\documentclass[11pt]{amsart}
\usepackage[T1]{fontenc}
\usepackage{lmodern}
\usepackage{amsmath,amssymb,amsthm,mathtools}
\usepackage[a4paper,margin=30mm]{geometry}
\usepackage[colorlinks=true,linkcolor=blue,citecolor=blue,urlcolor=blue]{hyperref}

\newtheorem{theorem}{Theorem}[section]
\newtheorem*{theorem*}{Theorem}
\newtheorem{lemma}[theorem]{Lemma}
\newtheorem{corollary}[theorem]{Corollary}
\newtheorem{proposition}[theorem]{Proposition}
\newtheorem*{proposition*}{Proposition}
\newtheorem{remark}[theorem]{Remark}

\theoremstyle{remark}
\newtheorem{fact}[theorem]{Fact}
\title[$\mathcal{PGF}=\mathcal{GP}$ for algebras finite over commutative Noetherian rings]{$\mathcal{PGF}=\mathcal{GP}$ for algebras finite over commutative Noetherian rings}
\author{Yuhui Hua\textsuperscript{*}}
\address{School of Mathematical Sciences, Shanghai Jiao Tong University, Shanghai 200240, China}
\email{huayuhui@sjtu.edu.cn}
\author{Shun-Jie Li}
\address{School of Mathematical Sciences, Shanghai Jiao Tong University, Shanghai 200240, China}
\email{lishunjie@sjtu.edu.cn}
\author{Chencheng Zhang}
\thanks{\textsuperscript{*}Corresponding author.}
\address{School of Mathematical Sciences, Shanghai Jiao Tong University, Shanghai 200240, China}
\email{zhangchencheng@sjtu.edu.cn}
\thanks{The authors were supported by the National Natural Science Foundation of China (No.~12131015).}

\begin{document}
\begin{abstract}
Let $S$ be a commutative Noetherian ring, $A$ be a finite $S$-algebra, i.e., $A$ is finitely generated as an $S$-module. 
We prove that an $A$-module is projectively coresolved Gorenstein flat if and only if it is Gorenstein projective. In other words, $\mathcal{PGF}(A)=\mathcal{GP}(A)$.
\end{abstract}
\maketitle

\section{Introduction}

Let $R$ be a unital ring.
A module ${}_RM$ is \emph{Gorenstein projective} if it is a cycle of an acyclic projective complex $\mathbf P$ such that $\mathrm{Hom}_R(\mathbf P,Q)$ is acyclic for each projective module ${}_RQ$.
It is \emph{Gorenstein flat} if it is a cycle of an acyclic flat complex $\mathbf F$ such that $J\otimes_R\mathbf F$ is acyclic for each injective module $J_R$.
It is \emph{projectively coresolved Gorenstein flat} if it is a cycle of an acyclic projective complex $\mathbf P$ such that $J\otimes_R\mathbf P$ is acyclic for each injective module $J_R$.
We denote the classes of Gorenstein projective, Gorenstein flat, and projectively coresolved Gorenstein flat modules by $\mathcal{GP}(R)$, $\mathcal{GF}(R)$, and $\mathcal{PGF}(R)$, respectively. Write $\mathrm{Proj}(R)$ and $\mathrm{Flat}(R)$ for the classes of projective and flat left $R$-modules, respectively.

\medskip

Gorenstein projective modules extend Auslander--Bridger's modules of $G$-dimension zero \cite{AB}, see \cite{EJ} for a modern definition.
Enochs, Jenda, and Torrecillas introduced Gorenstein flat modules in 1993 \cite{EJT}.
It is known that $(\mathcal{GP}(R),\mathcal{GP}(R)^{\perp_1})$ is a hereditary cotorsion pair \cite[Corollary~3.4(1)]{cortes-saroch};
whether $(\mathcal{GP}(R),\mathcal{GP}(R)^{\perp_1})$ is complete, or whether $\mathcal{GP}(R)$ is precovering, remains open in general \cite[Question~7.2]{cox}.
In 2020, \v{S}aroch and \v{S}\v{t}ov\'{i}\v{c}ek introduced projectively coresolved Gorenstein flat modules,
and proved that $(\mathcal{PGF}(R),\mathcal{PGF}(R)^{\perp_1})$ is a complete hereditary cotorsion pair \cite[Theorem~4.9]{saroch-stovicek}, and thus $\mathcal{PGF}(R)$ is special precovering.
Since $\mathcal{PGF}(R)\subseteq\mathcal{GP}(R)$ and $\mathcal{PGF}(R)$ enjoys good homological properties, it is natural to ask whether $\mathcal{PGF}(R)=\mathcal{GP}(R)$.
By \cite{iacob,saroch-stovicek}, this is equivalent to Holm's earlier question whether $\mathcal{GP}(R)\subseteq\mathcal{GF}(R)$. See \cite{holm-thesis} (also recorded in \cite[Remark~4.5(4)]{ding-li-mao}).
Recently, a Gorenstein projective module which is not Gorenstein flat, is found over a commutative coherent ring of Krull dimension zero \cite{zhang}.

\medskip

Several sufficient conditions for $\mathcal{GP}(R)\subseteq\mathcal{GF}(R)$ have been given.
In 2004, Holm proved $\mathcal{GP}(R)\subseteq\mathcal{GF}(R)$ for right coherent rings of finite left big finitistic projective dimension $\operatorname{FPD}(R)$ \cite[Proposition~3.4]{holm}.
In 2011, Murfet and Salarian proved the equivalence of total acyclicity and $F$-total acyclicity for complexes of projective modules over commutative Noetherian rings of finite Krull dimension \cite[Lemma~4.20]{murfet-salarian}.
In 2018, Christensen and Kato proved that totally acyclic complexes of projective modules are $F$-totally acyclic for commutative coherent rings $R$ with $\sup\{\operatorname{pd}_R F\mid F\text{ is flat}\}<\infty$. See \cite[Claim~2.5]{christensen-kato}.
They also proved it for commutative Noetherian rings $R$ with $\operatorname{pd}_R R_{\mathfrak p}<\infty$ for each $\mathfrak p\in\operatorname{Spec}R$. See \cite[Claim~2.8]{christensen-kato}.
In 2020, Iacob proved that $\mathcal{PGF}(R)=\mathcal{GP}(R)$ over any ring $R$ for which each injective right $R$-module has finite flat dimension \cite{iacob}.
Also in 2020, \v{S}aroch and \v{S}\v{t}ov\'{i}\v{c}ek showed that, over a right $\aleph_0$-coherent ring $R$, $\mathcal{PGF}(R^{\mathrm{op}})$ are precisely the modules filtered by countably presented Gorenstein projective right modules \cite[Theorem~4.9 and the remark following Example~4.10]{saroch-stovicek}.

\medskip

In this paper, let $S$ be a commutative Noetherian ring, $A$ be a finite $S$-algebra, i.e., $A$ is finitely generated as an $S$-module. 
We prove that an $A$-module is projectively coresolved Gorenstein flat if and only if it is Gorenstein projective. In other words, $\mathcal{PGF}(A)=\mathcal{GP}(A)$.

\begin{theorem*}
Let $S$ be a commutative Noetherian ring and $A$ be an $S$-algebra which is finitely generated as an $S$-module. Let $\mathbf P$ be an acyclic complex of projective left $A$-modules. Then the following are equivalent.

$(1)$ \ $\mathrm{Hom}_A(\mathbf P,Q)$ is acyclic for each $Q\in\mathrm{Proj}(A)$.

$(2)$ \ $J\otimes_A\mathbf P$ is acyclic for each injective right $A$-module $J$.

In particular, $\mathcal{PGF}(A)=\mathcal{GP}(A)$.
\end{theorem*}

It is known that $(2)$ implies $(1)$ \cite[Theorem~4.4]{saroch-stovicek}. We now prove the converse.

\section{Proof of the main theorem}

Fix $\mathbf P$ such that $\mathrm{Hom}_S(\mathbf P, Q)$ is acyclic for each $Q \in \mathrm{Proj}(A)$. Since $A$ is a Noetherian $S$-module, $A$ is left and right Noetherian.

\begin{lemma}\label{lem_Ra}
For $a\in S$, write $S_a=\{1,a,a^2,\ldots\}^{-1}S$. Then
$\mathrm{Hom}_A(\mathbf P,S_a\otimes_S A)$ is acyclic.
\end{lemma}

\begin{proof}
Apply $\mathrm{Hom}_A(\mathbf P,-)$ to the exact sequence
\begin{equation*}
0\longrightarrow A^{(\aleph_0)}
\xrightarrow{(e_i-ae_{i+1})} A^{(\aleph_0)}
\longrightarrow S_a\otimes_S A\longrightarrow 0
\end{equation*}
and use the fundamental theorem of homological algebra.
\end{proof}

Choose an arbitrary maximal ideal of $S$ together with a set of generators, writing $\mathfrak m = (a_1, a_2,\ldots, a_n)$. Write the \v{C}ech complex as $\check C _\mathfrak m = \bigotimes_{i=1}^n [S \xrightarrow{\text{canonical}} S_{a_i}]$.

\begin{lemma}\label{lem_Cm}
$\mathrm{Hom}_A(\mathbf P,\check C_{\mathfrak m}\otimes_S A)$ is acyclic.
\end{lemma}

\begin{proof}
Note that $S_a \otimes _S S_b \cong S_{ab}$. By Lemma~\ref{lem_Ra}, $\mathrm{Hom}_S(\mathbf P, (\check C_{\mathfrak m}\otimes_S A)^{i})$ is acyclic for each $i\in \mathbb{Z}$. Since $\check C_{\mathfrak m}\otimes_S A$ is bounded, by induction, $\mathrm{Hom}_S(\mathbf P, \check C_{\mathfrak m})$ is acyclic.
\end{proof}

We now turn to $\widehat{S_\mathfrak m}$, the completion of the local ring $S_\mathfrak m$ at its maximal ideal. Write
\begin{equation*}
(-)_{\widehat{S_\mathfrak m}} : \mathrm{Mod}_{S} \longrightarrow \mathrm{Mod}_{\widehat{S_\mathfrak m}}, \qquad X \longmapsto \widehat{S_\mathfrak m} \otimes _S X.
\end{equation*}
Set $A_{\widehat{S_{\mathfrak m}}}=\widehat{S_{\mathfrak m}}\otimes_S A$.
It is a $\widehat{S_{\mathfrak m}}$-algebra which is finitely generated as a $\widehat{S_{\mathfrak m}}$-module and hence is
left and right Noetherian. The above functor induces
\begin{equation*}
(-)_{\widehat{S_{\mathfrak m}}}\colon
\mathrm{Mod}_A\longrightarrow
\mathrm{Mod}_{A_{\widehat{S_{\mathfrak m}}}},
\qquad X\longmapsto\widehat{S_{\mathfrak m}}\otimes_S X.
\end{equation*}

\begin{fact}\label{fact}
If each element of $X \in \mathrm{Mod}_S$ is killed by some power of $\mathfrak m$, then the canonical map $X \to X_{\widehat{S_\mathfrak m}}$ is an isomorphism \cite[Lemme~1]{BeauvilleLaszlo}. The functor $(-)_{\widehat{S_\mathfrak m}}$ is exact, that is, $\widehat{S_\mathfrak m}$ is a flat $S$-module \cite[Proposition~10.14]{AtiyahMacdonald}.
\end{fact}

\begin{lemma}\label{lem_vCm}
$\mathrm{Hom}_{A_{\widehat{S_{\mathfrak m}}}}(\mathbf P_{\widehat{S_{\mathfrak m}}},(\check C_{\mathfrak m}\otimes_S A)_{\widehat{S_{\mathfrak m}}})$ is acyclic.
\end{lemma}

\begin{proof}
By \cite[Theorem~5.1.20]{brodmann-sharp}, each element of the cohomology of $\check C_{\mathfrak m}\otimes_S A$ is killed by some power of $\mathfrak m$. By Fact~\ref{fact}, the canonical map $\check C_{\mathfrak m}\otimes_S A\longrightarrow(\check C_{\mathfrak m}\otimes_S A)_{\widehat{S_{\mathfrak m}}}$ is a quasi-isomorphism. Since $\mathbf P$ is a projective complex, the functor $\mathrm{Hom}_A(\mathbf P,-)$ preserves quasi-isomorphisms of bounded complexes. Then by Lemma~\ref{lem_Cm}, $\mathrm{Hom}_A(\mathbf P,(\check C_{\mathfrak m}\otimes_S A)_{\widehat{S_{\mathfrak m}}})$ is acyclic. Finally, we have
\begin{align*}
\mathrm{Hom}_A(\mathbf P,(\check C_{\mathfrak m}\otimes_S A)_{\widehat{S_{\mathfrak m}}})& \cong
\mathrm{Hom}_{A_{\widehat{S_{\mathfrak m}}}}
\bigl(A_{\widehat{S_{\mathfrak m}}}\otimes_A\mathbf P_{\widehat{S_{\mathfrak m}}},
(\check C_{\mathfrak m}\otimes_S A)_{\widehat{S_{\mathfrak m}}}\bigr)\\
&\cong\mathrm{Hom}_{A_{\widehat{S_{\mathfrak m}}}}
\bigl(\mathbf P_{\widehat{S_{\mathfrak m}}},
(\check C_{\mathfrak m}\otimes_S A)_{\widehat{S_{\mathfrak m}}}\bigr).\qedhere  
\end{align*}  
\end{proof}

In the following discussion, we identify $(\check C_{\mathfrak m})_{\widehat{S_\mathfrak m}}$ with $\check C_{\mathfrak m \widehat{S_\mathfrak m}}$ by the following isomorphisms
\begin{equation*}
(\check C_{\mathfrak m})_{\widehat{S_\mathfrak m}} = \widehat{S_\mathfrak m}\otimes _S \check C_{\mathfrak m}  \cong \bigotimes_{i=1}^n [\widehat{S_\mathfrak m} \xrightarrow{\text{canonical}} (\widehat{S_\mathfrak m})_{a_i}] = \check C_{\mathfrak m \widehat{S_\mathfrak m}}.
\end{equation*}
Here $\widehat{S_\mathfrak m}$ is a local ring with maximal ideal $\mathfrak{m}\widehat{S_\mathfrak m}$, and $\check C_{\mathfrak m \widehat{S_\mathfrak m}}$ is the \v Cech complex induced by $(a_1,a_2,\ldots,a_n)$. Let $\kappa$ be the residue field of $\widehat{S_\mathfrak m}$.
Write $E = E_{\widehat{S_\mathfrak m}}(\kappa)$ and the duality
\begin{equation*}
(-)^\vee\colon
\bigl(\mathrm{Mod}_{A_{\widehat{S_{\mathfrak m}}}}\bigr)^{\mathrm{op}}
\longrightarrow
\mathrm{Mod}_{(A_{\widehat{S_{\mathfrak m}}})^{\mathrm{op}}},
\qquad
X\longmapsto \mathrm{Hom}_{\widehat{S_{\mathfrak m}}}(X,E).
\end{equation*}
Here $E$ is an injective cogenerator \cite[p.~528]{matlis}, so $(-)^\vee$ preserves and reflects exactness.

\begin{lemma}\label{lem_dual_acyclic}
$(\check C_{\mathfrak m\widehat{S_{\mathfrak m}}}\otimes_S A)^\vee
\otimes_{A_{\widehat{S_{\mathfrak m}}}}\mathbf P_{\widehat{S_{\mathfrak m}}}$ is acyclic.
\end{lemma}

\begin{proof}
First we show that the canonical map $\check C_{\mathfrak m\widehat{S_{\mathfrak m}}}\otimes_S A \longrightarrow(\check C_{\mathfrak m\widehat{S_{\mathfrak m}}}\otimes_S A)^{\vee\vee}$ is a quasi-isomorphism.
By \cite[Theorem~7.1.3]{brodmann-sharp}, the cohomology groups of $\check C_{\mathfrak m\widehat{S_{\mathfrak m}}}\otimes_S A$ are Artinian, since they are the local cohomology groups of the finitely generated $\widehat{S_{\mathfrak m}}$-module $A_{\widehat{S_{\mathfrak m}}}$ at the maximal ideal $\mathfrak m\widehat{S_{\mathfrak m}}$. 
Moreover, since $\widehat{S_\mathfrak m}$ is complete Noetherian, by \cite[Corollary~4.3(3)]{matlis}, we have $H^i ( \check C_{\mathfrak m \widehat{S_\mathfrak m}}\otimes_S A) \cong (H^i( \check C_{\mathfrak m \widehat{S_\mathfrak m}}\otimes_S A))^{\vee\vee}$ for each $i\in\mathbb{Z}$. 
Since $(-)^{\vee}$ is exact, we have $(H^i( \check C_{\mathfrak m \widehat{S_\mathfrak m}}\otimes_S A))^{\vee\vee} \cong H^i(( \check C_{\mathfrak m \widehat{S_\mathfrak m}}\otimes_S A)^{\vee\vee})$ for each $i \in\mathbb{Z}$.
It follows that 
\[\mathrm{Hom}_{\widehat{S_\mathfrak m}}(\mathbf P_{\widehat{S_\mathfrak m}},  \check C_{\mathfrak m\widehat{S_{\mathfrak m}}}\otimes_S A)\longrightarrow
\mathrm{Hom}_{\widehat{S_\mathfrak m}}(\mathbf P_{\widehat{S_\mathfrak m}}, (\check C_{\mathfrak m\widehat{S_{\mathfrak m}}}\otimes_S A)^{\vee\vee})\]
is a quasi-isomorphism.
By Lemma~\ref{lem_vCm}, \(\mathrm{Hom}_{A_{\widehat{S_{\mathfrak m}}}}(\mathbf P_{\widehat{S_{\mathfrak m}}},(\check C_{\mathfrak m\widehat{S_{\mathfrak m}}}\otimes_S A)^{\vee\vee})\)
is an acyclic complex.
By Tensor--Hom adjunction,
\begin{equation*}
 \bigl((\check C_{\mathfrak m\widehat{S_{\mathfrak m}}}\otimes_S A)^\vee
\otimes_{A_{\widehat{S_{\mathfrak m}}}}\mathbf P_{\widehat{S_{\mathfrak m}}}\bigr)^\vee \cong \mathrm{Hom}_{A_{\widehat{S_{\mathfrak m}}}}(\mathbf P_{\widehat{S_{\mathfrak m}}},(\check C_{\mathfrak m\widehat{S_{\mathfrak m}}}\otimes_S A)^{\vee\vee}).
\end{equation*}
Since $(-)^\vee$ reflects exactness, $(\check C_{\mathfrak m\widehat{S_{\mathfrak m}}}\otimes_S A)^\vee
\otimes_{A_{\widehat{S_{\mathfrak m}}}}\mathbf P_{\widehat{S_{\mathfrak m}}}$ is acyclic.
\end{proof}

\begin{proposition}\label{prop_key}
$J\otimes_{A_{\widehat{S_{\mathfrak m}}}}
\mathbf P_{\widehat{S_{\mathfrak m}}}$ is acyclic for each injective
right $A_{\widehat{S_{\mathfrak m}}}$-module $J$.
\end{proposition}

\begin{proof}
\noindent\textbf{Step 1.}
We claim that such $J$ is a summand of 
\[\coprod\limits_{\mathfrak{q}\in \mathrm{Spec}(\widehat{S_{\mathfrak m}})} \left(\mathrm{Hom}_{\widehat{S_{\mathfrak m}}}
\bigl(A_{\widehat{S_{\mathfrak m}}},
E_{\widehat{S_{\mathfrak m}}}(\widehat{S_{\mathfrak m}}/\mathfrak q)\bigr)\right)^{(\kappa_{\mathfrak{q}})}.\]
Indeed, we regard $J$ as an $\widehat{S_{\mathfrak m}}$-module and take its injective envelope $J\hookrightarrow K$. 
Consider
\[\mathrm{Hom}_{\widehat{S_{\mathfrak m}}}(J,K)\cong \mathrm{Hom}_{A_{\widehat{S_{\mathfrak m}}}}(J,\mathrm{Hom}_{\widehat{S_{\mathfrak m}}}
(A_{\widehat{S_{\mathfrak m}}},K)),\]
the monomorphism $J\hookrightarrow K$ corresponds to a monomorphism $J\hookrightarrow \mathrm{Hom}_{\widehat{S_{\mathfrak m}}}
(A_{\widehat{S_{\mathfrak m}}},K)$, which splits since $J$ is an injective $A_{\widehat{S_{\mathfrak m}}}$-module.
By Matlis decomposition, $K$ is a direct sum of copies of $E_{\widehat{S_{\mathfrak m}}}(\widehat{S_{\mathfrak m}}/\mathfrak q)$'s.
Since $A_{\widehat{S_{\mathfrak m}}}$ is finitely generated over
$\widehat{S_{\mathfrak m}}$, it follows that $\mathrm{Hom}_{\widehat{S_{\mathfrak m}}}(A_{\widehat{S_{\mathfrak m}}},-)$ commutes with direct sums. 

By claim, in the rest of the proof, we may assume $J=\mathrm{Hom}_{\widehat{S_{\mathfrak m}}}
\bigl(A_{\widehat{S_{\mathfrak m}}},
E_{\widehat{S_{\mathfrak m}}}(\widehat{S_{\mathfrak m}}/\mathfrak q)\bigr)$ for some $\mathfrak{q}\in \mathrm{Spec}(\widehat{S_{\mathfrak m}})$.
\medskip

\noindent\textbf{Step 2.} 
By Lemma~\ref{lem_dual_acyclic}, $(\check C_{\mathfrak m\widehat{S_{\mathfrak m}}}\otimes_S A)^\vee
\otimes_{A_{\widehat{S_{\mathfrak m}}}}\mathbf P_{\widehat{S_{\mathfrak m}}}$ is acyclic. 
After localization, we have an isomorphism between acyclic complexes
\begin{equation*}
\begin{aligned}
\Bigl((\check C_{\mathfrak m\widehat{S_{\mathfrak m}}}\otimes_S A)^\vee \otimes_{A_{\widehat{S_{\mathfrak m}}}}
\mathbf P_{\widehat{S_{\mathfrak m}}}
\Bigr)_{\mathfrak q}
\cong
\bigl((\check C_{\mathfrak m\widehat{S_{\mathfrak m}}}\otimes_S A)^\vee\bigr)_{\mathfrak q}
\otimes_{(A_{\widehat{S_{\mathfrak m}}})_{\mathfrak q}}
(\mathbf P_{\widehat{S_{\mathfrak m}}})_{\mathfrak q}. 
\end{aligned}
\end{equation*}
Since $\check C_{\mathfrak q(\widehat{S_{\mathfrak m}})_{\mathfrak q}}$ is a bounded complex of flat $(\widehat{S_{\mathfrak m}})_{\mathfrak q}$-modules,
applying $\check C_{\mathfrak q(\widehat{S_{\mathfrak m}})_{\mathfrak q}}\otimes_{(\widehat{S_{\mathfrak m}})_{\mathfrak q}}-$,
we obtain an acyclic complex \(\check C_{\mathfrak q(\widehat{S_{\mathfrak m}})_{\mathfrak q}}
\otimes_{(\widehat{S_{\mathfrak m}})_{\mathfrak q}}\Bigl(\bigl((\check C_{\mathfrak m\widehat{S_{\mathfrak m}}}\otimes_S A)^\vee\bigr)_{\mathfrak q}
\otimes_{(A_{\widehat{S_{\mathfrak m}}})_{\mathfrak q}}(\mathbf P_{\widehat{S_{\mathfrak m}}})_{\mathfrak q}\Bigr).\)

\medskip

\noindent\textbf{Step 3.}
We want to replace $\check C_{\mathfrak q(\widehat{S_{\mathfrak m}})_{\mathfrak q}}
\otimes_{(\widehat{S_{\mathfrak m}})_{\mathfrak q}}\bigl((\check C_{\mathfrak m\widehat{S_{\mathfrak m}}}\otimes_S  A)^\vee\bigr)_{\mathfrak q}$ with $J$.
By \cite[(2.4.2)]{lipman}, $(\check C_{\mathfrak m\widehat{S_{\mathfrak m}}})^\vee$ is a normalized dualizing complex over $\widehat{S_{\mathfrak m}}$.
Then by \cite[Remark~2.16(1)]{hrbek-nakamura-stovicek}, $\bigl((\check C_{\mathfrak m\widehat{S_{\mathfrak m}}})^\vee\bigr)_{\mathfrak q}$ is a dualizing complex over $(\widehat{S_{\mathfrak m}})_{\mathfrak q}$.
Write $\bigl((\check C_{\mathfrak m\widehat{S_{\mathfrak m}}})^\vee\bigr)_{\mathfrak q}
=\omega[s]$, a shift of the normalized dualizing complex $\omega$.
By \cite[Proposition~3.1.2]{lipman}, there is an isomorphism of derived functors
\[
\check C_{\mathfrak q(\widehat{S_{\mathfrak m}})_{\mathfrak q}}
\otimes_{(\widehat{S_{\mathfrak m}})_{\mathfrak q}}(-)
\cong
\mathbf R\Gamma_{\mathfrak q(\widehat{S_{\mathfrak m}})_{\mathfrak q}}(-),
\]
and \cite[\S~2.4]{lipman} gives
\[
\mathbf R\Gamma_{\mathfrak q(\widehat{S_{\mathfrak m}})_{\mathfrak q}}(\omega)
\cong
E_{(\widehat{S_{\mathfrak m}})_{\mathfrak q}}
\bigl((\widehat{S_{\mathfrak m}})_{\mathfrak q}/
\mathfrak q(\widehat{S_{\mathfrak m}})_{\mathfrak q}\bigr)[0].
\]
By \cite[Theorem~3.6]{matlis}, there is an isomorphism of
$\widehat{S_{\mathfrak m}}$-modules
\[E_{(\widehat{S_{\mathfrak m}})_{\mathfrak q}}
\bigl((\widehat{S_{\mathfrak m}})_{\mathfrak q}/
\mathfrak q(\widehat{S_{\mathfrak m}})_{\mathfrak q}\bigr)
\cong
E_{\widehat{S_{\mathfrak m}}}(\widehat{S_{\mathfrak m}}/\mathfrak q).\]
Therefore we have an isomorphism in $\mathrm D((\widehat{S_{\mathfrak m}})_{\mathfrak q})$
\begin{equation}\label{eq_homot_equiv}
\check C_{\mathfrak q(\widehat{S_{\mathfrak m}})_{\mathfrak q}}
\otimes_{(\widehat{S_{\mathfrak m}})_{\mathfrak q}}
\bigl((\check C_{\mathfrak m\widehat{S_{\mathfrak m}}})^\vee\bigr)_{\mathfrak q}
\cong
E_{\widehat{S_{\mathfrak m}}}(\widehat{S_{\mathfrak m}}/\mathfrak q)[s].
\end{equation}
Since both sides of \eqref{eq_homot_equiv} are bounded complexes of injective modules, \eqref{eq_homot_equiv} is represented by a homotopy equivalence.
Applying $\mathrm{Hom}_{(\widehat{S_{\mathfrak m}})_{\mathfrak q}} ((A_{\widehat{S_{\mathfrak m}}})_{\mathfrak q},-)$ to \eqref{eq_homot_equiv}, the right side is $J[s]$ by Tensor--Hom adjunction. For the left side, we first use the Tensor--Hom adjunction
\[(\check C_{\mathfrak m\widehat{S_{\mathfrak m}}}\otimes_S A)^\vee
\cong\mathrm{Hom}_{\widehat{S_{\mathfrak m}}}
\bigl(A_{\widehat{S_{\mathfrak m}}},
      (\check C_{\mathfrak m\widehat{S_{\mathfrak m}}})^\vee\bigr).\]
Since $A_{\widehat{S_{\mathfrak m}}}$ is finitely presented over $\widehat{S_{\mathfrak m}}$, one has the isomorphisms after localizition at $\mathfrak q$
\[
\bigl((\check C_{\mathfrak m\widehat{S_{\mathfrak m}}}
         \otimes_S A)^\vee\bigr)_{\mathfrak q}
\cong
\left(\mathrm{Hom}_{\widehat{S_{\mathfrak m}}}
\bigl(A_{\widehat{S_{\mathfrak m}}},
      (\check C_{\mathfrak m\widehat{S_{\mathfrak m}}})^\vee\bigr)\right)_{\mathfrak{q}}
\cong
\mathrm{Hom}_{(\widehat{S_{\mathfrak m}})_{\mathfrak q}}
\Bigl(
 (A_{\widehat{S_{\mathfrak m}}})_{\mathfrak q},
 \bigl((\check C_{\mathfrak m\widehat{S_{\mathfrak m}}})^\vee
       \bigr)_{\mathfrak q}
\Bigr).\]
Since $\check C_{\mathfrak q(\widehat{S_{\mathfrak m}})_{\mathfrak q}}$ is a bounded complex of flat $(\widehat{S_{\mathfrak m}})_{\mathfrak q}$-modules and $(A_{\widehat{S_{\mathfrak m}}})_{\mathfrak q}$ is finitely presented over $(\widehat{S_{\mathfrak m}})_{\mathfrak q}$, we have an isomorphism of complexes
\[\begin{aligned}
&\check C_{\mathfrak q(\widehat{S_{\mathfrak m}})_{\mathfrak q}}
\otimes_{(\widehat{S_{\mathfrak m}})_{\mathfrak q}}
\mathrm{Hom}_{(\widehat{S_{\mathfrak m}})_{\mathfrak q}}
\Bigl(
 (A_{\widehat{S_{\mathfrak m}}})_{\mathfrak q},
 \bigl((\check C_{\mathfrak m\widehat{S_{\mathfrak m}}})^\vee
       \bigr)_{\mathfrak q}
\Bigr)\\
&\qquad\cong
\mathrm{Hom}_{(\widehat{S_{\mathfrak m}})_{\mathfrak q}}
\Bigl(
 (A_{\widehat{S_{\mathfrak m}}})_{\mathfrak q},
 \check C_{\mathfrak q(\widehat{S_{\mathfrak m}})_{\mathfrak q}}
 \otimes_{(\widehat{S_{\mathfrak m}})_{\mathfrak q}}
 \bigl((\check C_{\mathfrak m\widehat{S_{\mathfrak m}}})^\vee
       \bigr)_{\mathfrak q}
\Bigr).
\end{aligned}\]
Therefore we obtain a homotopy equivalence of complexes of right
$(A_{\widehat{S_{\mathfrak m}}})_{\mathfrak q}$-modules
\[
\check C_{\mathfrak q(\widehat{S_{\mathfrak m}})_{\mathfrak q}}
\otimes_{(\widehat{S_{\mathfrak m}})_{\mathfrak q}}
\bigl((\check C_{\mathfrak m\widehat{S_{\mathfrak m}}}
       \otimes_S A)^\vee\bigr)_{\mathfrak q}
\cong J[s].
\]

\medskip

\noindent\textbf{Step 4.}
We know from \textbf{Step~2} that  
$\check C_{\mathfrak q(\widehat{S_{\mathfrak m}})_{\mathfrak q}}
\otimes_{(\widehat{S_{\mathfrak m}})_{\mathfrak q}} \bigl((\check C_{\mathfrak m\widehat{S_{\mathfrak m}}}
\otimes_S A)^\vee\bigr)_{\mathfrak q}
\otimes_{A_{\widehat{S_{\mathfrak m}}}}(\mathbf P_{\widehat{S_{\mathfrak m}}})_{\mathfrak q}$
is acyclic, hence by \textbf{Step~3},
$J\otimes_{(A_{\widehat{S_{\mathfrak m}}})_{\mathfrak q}}
(\mathbf P_{\widehat{S_{\mathfrak m}}})_{\mathfrak q}$ is acyclic.
Consider
\begin{equation*}
J\otimes_{(A_{\widehat{S_{\mathfrak m}}})_{\mathfrak q}}
(\mathbf P_{\widehat{S_{\mathfrak m}}})_{\mathfrak q}
\cong
J\otimes_{(A_{\widehat{S_{\mathfrak m}}})_{\mathfrak q}}
\bigl((A_{\widehat{S_{\mathfrak m}}})_{\mathfrak q}
\otimes_{A_{\widehat{S_{\mathfrak m}}}}
\mathbf P_{\widehat{S_{\mathfrak m}}}\bigr)
\cong
J\otimes_{A_{\widehat{S_{\mathfrak m}}}}
\mathbf P_{\widehat{S_{\mathfrak m}}}.
\end{equation*}
Hence $J\otimes_{A_{\widehat{S_{\mathfrak m}}}}
\mathbf P_{\widehat{S_{\mathfrak m}}}$ is acyclic.  Proposition~\ref{prop_key} now follows from the reduction in \textbf{Step~1}. \qedhere 
\end{proof}

\begin{remark}
The proposition shows that $\mathcal{PGF}(A_{\widehat{S_{\mathfrak m}}})=\mathcal{GP}(A_{\widehat{S_{\mathfrak m}}})$.
\end{remark}

Write $A_{\mathfrak m}=S_{\mathfrak m}\otimes_S A$ and $\mathbf P_{\mathfrak m}=S_{\mathfrak m}\otimes_S\mathbf P$.

\begin{lemma}\label{lem_PRm}
$J\otimes_{A_{\mathfrak m}}\mathbf P_{\mathfrak m}$ is acyclic for each injective right $A_{\mathfrak m}$-module $J$.
\end{lemma}

\begin{proof}
Since $\widehat{S_{\mathfrak m}}$ is faithfully flat over
$S_{\mathfrak m}$, the canonical map $J\longrightarrow\widehat{S_{\mathfrak m}}\otimes_{S_{\mathfrak m}}J$ is injective. 
Take the injective envelope $\widehat{S_{\mathfrak m}}\otimes_{S_{\mathfrak m}}J\longrightarrow K$ over $A_{\widehat{S_{\mathfrak m}}}$. 
Then $J$ is a direct summand of $K$ as a right $A_{\mathfrak m}$-module. 
Thus $J\otimes_{A_{\mathfrak m}}\mathbf P_{\mathfrak m}$ is a direct summand of \(
K\otimes_{A_{\mathfrak m}}\mathbf P_{\mathfrak m}\cong K\otimes_{A_{\widehat{S_{\mathfrak m}}}}\mathbf P_{\widehat{S_{\mathfrak m}}},
\)
which is acyclic by Proposition~\ref{prop_key}.
\end{proof}

\begin{lemma}\label{lem_global}
$J\otimes_A\mathbf P$ is acyclic for each injective right $A$-module $J$.
\end{lemma}

\begin{proof}
Since $A$ is right Noetherian, $J_{\mathfrak m}$ is an injective right $A_{\mathfrak m}$-module by \cite[Lemma~1.2(b)]{bass}.
By Lemma~\ref{lem_PRm}, $(J\otimes_A\mathbf P)_{\mathfrak m}\cong J_{\mathfrak m}\otimes_{A_{\mathfrak m}}\mathbf P_{\mathfrak m}$ is acyclic for each maximal ideal $\mathfrak m$ of $S$, and hence $J\otimes_A\mathbf P$ is acyclic.
\end{proof}

\begin{corollary}
Let $S$ be a commutative Noetherian ring and $A$ be an $S$-algebra which is finitely generated as an $S$-module. Then the following are equivalent for each acyclic complex $\mathbf P$ of projective left $A$-modules.
\begin{enumerate}
\item $\mathrm{Hom}_A(\mathbf P,Q)$ is acyclic for each $Q\in\mathrm{Proj}(A)$.
\item $J\otimes_A\mathbf P$ is acyclic for each injective right $A$-module $J$.
\item $\mathrm{Hom}_A(\mathbf P,A^{(\aleph_0)})$ is acyclic.
\item $\mathrm{Hom}_A(\mathbf P,F)$ is acyclic for each $F\in\mathrm{Flat}(A)$.
\end{enumerate}
\end{corollary}

\begin{proof}
The implications $(4)\Rightarrow(1)\Rightarrow(3)$ are immediate.
Under $(3)$, Lemma~\ref{lem_Ra} applies and the rest of the proof of the main theorem gives $(2)$.
Under $(2)$, each $Z^i\mathbf P$ is a $\mathcal{PGF}$ module. 
By \cite[Theorem~4.4]{saroch-stovicek}, $\mathrm{Ext}^1_A(Z^i\mathbf P,F)=0$ for each flat left $A$-module $F$ and each $i\in\mathbb Z$. 
Applying $\mathrm{Hom}_A(-,F)$ to each $0\longrightarrow Z^i\mathbf P\longrightarrow P^i \longrightarrow Z^{i+1}\mathbf P\longrightarrow 0$ gives $(4)$.
\end{proof}

\end{document}